\documentclass[a4paper]{article}

\usepackage{amssymb}
\usepackage{amsmath}

\newtheorem{theorem}{Theorem}[section]

\newtheorem{e-proposition}[theorem]{Proposition}

\newtheorem{e-definition}[theorem]{Definition\rm}

\newtheorem{remarks}[theorem]{Remarks}

\def\og{\leavevmode\raise.3ex\hbox{$\scriptscriptstyle\langle\!\langle$~}}
\def\fg{\leavevmode\raise.3ex\hbox{~$\!\scriptscriptstyle\,\rangle\!\rangle$}}

\title{Addendum to "The $\hat A$-genus of $S^1$-manifolds with finite second homotopy group"}

\begin{document}


\author{Manuel Amann\\ Department of Mathematics, University of Augsburg, Germany \thanks{The first named author thanks the University of Augsburg for their continued support.}\and Anand Dessai\\ Department of Mathematics, University of Fribourg, Switzerland \thanks{The second named author was supported by the SNSF-Project 200020E\_193062 and the DFG-Priority programme SPP 2026.}}

\maketitle

\setcounter{section}{1}
Theorem 1.1 of \cite{AD} states that for any $k>1$ there exists a smooth simply connected closed $4k$-dimen\-sio\-nal manifold with finite $\pi_2$ and non-vanishing $\hat A$-genus which admits a smooth effective $S^1$-action. The proof of this theorem was based on a general surgery lemma  \cite[Lemma 2.1]{AD}. As has been pointed out to us by Michael Wiemeler, the surgery lemma is not true in the form stated. Additional assumptions are necessary to ensure that the surgery can be done equivariantly. The purpose of this note is to remedy this situation and to give a proof of the theorem above using explicit equivariant surgeries. For $k> 2$ the proof is essentially the one given in \cite{AD}, for $k=2$ we apply \cite[Thm. 1.3]{W}.

\begin{theorem}\label{thm 1.1} For any $k>1$ there exists a smooth simply connected $4k$-dimen\-sio\-nal $\pi_2$-finite manifold $M^\prime$ with smooth effective $S^1$-action and $\hat A(M^\prime)\neq 0$.\end{theorem}

\noindent
{\bf Proof:} As in \cite{AD} we consider $M:={\mathbb C} P^{2k}$ with a linear $S^1$-action. More precisely, we consider the $S^1$-action on $M$ induced by the linear action
\begin{align*} S^1\times {\mathbb C} ^{2k+1}&\to {\mathbb C} ^{2k+1}\\
(\lambda, (z_0,\ldots ,z_{2k-2},z_{2k-1},z_{2k}))&\mapsto (z_0,\ldots ,z_{2k-2}, \lambda \cdot z_{2k-1},\lambda ^{-1}\cdot z_{2k} ).\end{align*}
The connected components of the fixed point set $M^{S^1}$ are two isolated fixed points and $N:={\mathbb C} P^{2k-2}$.

We first assume that the dimension of $N$ is $>4$, i.e., that $k>2$. Let $\kappa :M\to BSO$ be a classifying map for the stable normal bundle of $M$. Consider a map $f:S^2\to N$ such that its composition with the inclusion map $N \hookrightarrow M$ represents a generator of the kernel of $\kappa _*:\pi _2(M)\to \pi _2(BSO)$. Since $\dim N >4$ we may assume that $f$ is a smooth embedding.

By construction, the normal bundle of $f(S^2) \subset M$, viewed as a real non-equivariant vector bundle, is trivial. It splits equivariantly as a direct sum $E\oplus F$, where $E$ is the normal bundle of $f(S^2)$ in $N$ and $F$ is the restriction to $f(S^2)$ of the normal bundle $\nu _N$  of $N$ in $M$. We note that $\nu _N$ can be identified with the direct sum $\overline \gamma \oplus \overline \gamma$, where $\gamma$ denotes the canonical complex line bundle over $N={\mathbb C} P^{2k-2}$. Under this identification, $\lambda \in S^1$ acts on one summand by multiplication with $\lambda $ and on the other summand by multiplication with $\lambda ^{-1}$. By changing the complex structure, we can identify $\nu _N$ equivariantly with the complex vector bundle $\gamma \oplus \overline \gamma$, where $\lambda \in S^1$ now acts by complex multiplication on both summands. In particular, $F$ can be identified equivariantly with the complex vector bundle $(\gamma \oplus \overline \gamma)\vert _{f(S^2)}$ which is trivial. Also, $E$ is a trivial real vector bundle over $f(S^2)$, since $E$ is stably trivial and has rank $\mathrm{rk}\, E=\dim N -2=4k-6>2$. Moreover, $S^1$ acts trivially on $E$. Hence, $F$ (resp. $E$) can be extended as a complex (resp. real) vector bundle over the disk $D^3$.

In conclusion, the normal bundle of $f(S^2)\subset M$ splits equivariantly as the direct sum of a trivial vector bundle $E$ with trivial $S^1$-action and a trivial vector bundle $F$ with non-trivial $S^1$-action, and the direct sum extends equivariantly to the disk $D^3$, where $S^1$ acts trivially on $D^3$. In particular, the sphere normal bundle of $f(S^2)\subset M$ can be extended equivariantly to the disk. Hence, one can perform $S^1$-equivariant surgery on $f(S^2)$. The result of the surgery is a simply connected $S^1$-manifold $M^\prime$ with $\pi _2(M^\prime)\cong {\mathbb Z}/2{\mathbb Z}$ which is $S^1$-equivariantly bordant to $M$.

It is well-known that the $\hat A$-genus does not vanish on ${\mathbb C} P^{2k}$. Since $M^\prime $ is bordant to ${\mathbb C} P^{2k}$ the same is true for $M^\prime $. This completes the argument if $k>2$.

For $k=2$,
one cannot argue as above, since the bundle $E$ does not extend to the disk. In fact, if $S^2\subset N={\mathbb C}  P ^2$ represents a non-trivial element of $\pi _2(M)$ then the Euler class of its normal bundle in $N$ is non-zero.

To construct $M^\prime$ for $k=2$ we will apply a different construction which involves non-equivariant surgery in the orbit space. This construction, which was pointed out to us by Michael Wiemeler, is quite general, the argument below is a special case of \cite[Thm. 1.3]{W}.

Let us consider $M={\mathbb C} P^{4}$ with the $S^1$-action as before. Let $M_{0}$ denote the union of principal orbits and $i:M_{0}\hookrightarrow M$ the inclusion map.  Note that the complement of $M_{0}$ is the union of $N$ and a 2-sphere which contains the two isolated fixed points. The restriction of the action to $M_{0}$ defines a principal $S^1$-bundle $p:M_{0} \to B_0$ with base space $B_0:=M_{0}/S^1$. Note that $i_*:\pi _2(M_{0})\to \pi _2(M)$ is an isomorphisms and $p_*:\pi _2(M_{0})\to \pi _2(B_0)$ is injective.

Consider a map $f:S^2\to M_{0}$ such that $i\circ f$ represents a generator of the kernel of $\kappa _*:\pi _2(M)\to \pi _2(BSO)$. Up to homotopy we may assume that $f$ and $\overline f:=p\circ f:S^2\to B_0$ are smooth embeddings. It follows that the restriction of the principal $S^1$-bundle to $\overline f(S^2)$ is trivial (with a section given by $f$) and that the normal bundle of $\overline f(S^2)\subset B_0$ is trivial. Hence, the disk normal bundle of $\overline f(S^2)$ can be identified with $\overline f(S^2)\times D^5$ and the disk normal bundle of $S^1\cdot f(S^2)=p^{-1}(\overline f(S^2))\cong S^1\times f(S^2)$ can be identified equivariantly with $S^1\times f(S^2)\times D^5$, where $S^1$ acts on $S^1$ by left multiplication and acts trivially on the other factors. Next perform surgery for $\overline f(S^2)\times D^5\subset B_0$ and make the corresponding modification for $S^1\times f(S^2)\times D^5\subset M_0$, i.e., consider
$$\Bigl(M_0 - (S^1\times f(S^2)\times \overset {\circ} {D^5})\Bigr) \cup \Bigl({S^1\times D^3\times S^4}\Bigr).$$ Let us denote the results of these modifications by $B_0^\prime $ and $M_0^\prime$, respectively, and let
$$M^\prime :=M_0^\prime\cup (M-M_0).$$ As before, the inclusion $M_{0}^\prime\hookrightarrow M^\prime$ induces an isomorphism on $\pi _2$ and the projection in the principal $S^1$-bundle $M_0^\prime \to B_0^\prime $ is injective on $\pi _2$. It follows that $\pi _2(M^\prime)\cong {\mathbb Z}/2{\mathbb Z}$. In addition, the long exact sequence of homotopy groups for the principal $S^1$-bundles together with the Seifert-van Kampen theorem shows that  $M^\prime$ is simply connected. By construction, $M^\prime $ is equivariantly bordant to $M={\mathbb C} P^{4}$. Hence, $\hat A(M^\prime)=\hat A({\mathbb C} P^{4})\neq 0$. This completes the argument for $k=2$.\hfill $\square$

\begin{remarks}
\begin{enumerate}
\item LeBrun and Salamon conjectured that any positive quaternionic K\"ahler manifold is symmetric. The conjecture has been proved by Hitchin, Poon-Salamon and LeBrun-Salamon in
dimensions $\leq 8$.  Hayde\'e and Rafael Herrera  \cite{HH} offered a proof in dimension $12$ involving a vanishing statement for the $\hat A$-genus of $\pi _2$-finite manifolds with $S^1$-action. The purpose of \cite[Thm. 1.1]{AD} was to show that this statement cannot be correct.
\item Recently, Buczy\'nski, Wi\'sniewski and Weber \cite{BWW} gave a proof of the LeBrun-Salamon conjecture in dimensions $12$ and $16$ by showing that the corresponding twistor spaces are adjoint varieties.
\item In \cite{W} Wiemeler uses equivariant surgery to construct examples of manifolds, as in Thm. \ref{thm 1.1}, which have non-trivial fundamental group. The examples have non-spin universal covering. On the other hand, Wiemeler shows that the elliptic genus is rigid for any closed even-dimensional $S^1$-manifold whose universal covering is spin. Using an argument of Hirzebruch and Slodowy \cite{HS} he concludes that the $\hat A$-genus of such a manifold vanishes if the action is non-trivial.
\end{enumerate}
\end{remarks}

\bigskip
We like to thank Michael Wiemeler for discussions and valuable comments.






\begin{thebibliography}{00}




\bibitem{AD} {M. Amann and A. Dessai}, {The $\hat A$-genus of $S^1$-manifolds with finite second homotopy group}, {C. R. Math. Acad. Sci. Paris} 348 (2010) 283--285
\bibitem{BWW} {J. Buczy\'nski, J. A. Wi\'sniewski, A. Weber}, {Algebraic torus actions on contact manifolds}, {J. Differential Geom.} 121 (2022) 227--289
\bibitem{HH} {H. and R. Herrera}, {$\hat A$-genus on non-spin manifolds with $S^1$ actions and the classification of positive quaternion-K\"ahler 12-manifolds}, J. Differential. Geom. 61 (2002) 341--364
\bibitem{HS} {F. Hirzebruch and P. Slodowy}, {Elliptic genera, involutions and homogeneous spin manifolds}, Geom. Ded. 35 (1990) 309--343
\bibitem{W} {M. Wiemeler}, {Rigidity of elliptic genera for non-spin manifolds}, preprint arXiv:2212.01059 (2022)
\end{thebibliography}
\end{document}